\documentclass[11pt]{article}
\usepackage{mathrsfs}
\usepackage{amssymb}
\usepackage{amsmath}
\usepackage{amsbsy}
\usepackage{epsfig}
\usepackage{enumerate}
\usepackage{bm}
\newtheorem{lemma}{Lemma}[section]
\newtheorem{theorem}{Theorem}[section]

\newtheorem{corollary}{Corollary}[section]
\newtheorem{definition}{Definition}[section]
\newtheorem{remark}{Remark}[section]

\newbox\TempBox \newbox\TempBoxA

\def\pr{\textsf{P}} 
\def\Sbep{\widehat{\mathbb E}} 
\def\cSbep{\widehat{\mathcal E}} 
\def\Capc{\mathbb V} 
\def\cCapc{\mathcal V} 

\def\underwiggle 1{
\ifmmode\setbox\TempBox=\hbox{$ 1$}\else\setbox\TempBox=\hbox{
1}\fi \setbox\TempBoxA=\hbox to \wd\TempBox{\hss\char'176\hss}
\rlap{\copy\TempBox}\smash{\lower9pt\hbox{\copy\TempBoxA}} }
\renewcommand{\baselinestretch}{1.7}

\begin{document}

\thispagestyle{empty}

\begin{center}
 { \LARGE\bf The convergence of the sums of independent random variables   under the sub-linear expectations $^{\ast}$}
\end{center}

\begin{center} {\sc
Li-Xin Zhang\footnote{This work was supported by grants from the NSF of China (Grant No. 11731012), the 973 Program (Grant No. 2015CB352302), Zhejiang Provincial Natural Science Foundation (Grant No. LY17A010016) and the Fundamental Research Funds for the Central Universities.
}
}\\
{\sl \small School of Mathematical Sciences, Zhejiang University, Hangzhou 310027} \\
(Email:stazlx@zju.edu.cn)\\
\end{center}

\renewcommand{\abstractname}{~}
\begin{abstract}
\centerline{\bf Abstract}
Let $\{X_n;n\ge 1\}$ be a sequence of independent random variables on a probability space $(\Omega, \mathcal{F}, \pr)$ and $S_n=\sum_{k=1}^n X_k$. It is well-known that the almost sure convergence, the convergence in probability and the convergence in distribution of $S_n$ are equivalent. In this paper, we prove   similar results for the independent random variables under the sub-linear expectations, and give a group of sufficient and necessary conditions for these convergence. For proving the results,   the Levy and Kolmogorov maximal inequalities for independent random variables under the sub-linear expectation are established. As an application of the maximal inequalities,  the sufficient and necessary conditions for the central limit theorem of independent and identically distributed random variables are also obtained.


\end{abstract}

\baselineskip 22pt

\renewcommand{\baselinestretch}{1.7}



\section{ Introduction and main results}\label{sect1}
\setcounter{equation}{0}

The convergence of the sums  of independent random variables are well-studied. For example, it is well-known that,  if $\{X_n;n\ge 1\}$ is a sequence of independent random variables on a probability space $(\Omega, \mathcal{F}, \pr)$, then that the infinite series     $\sum_{n=1}^{\infty}X_n$ is convergent almost surely, that it is convergent in probability and that it is convergent in distribution are equivalent. In this paper, we consider this elementary equivalence under the sub-linear expectations.
The general framework of the sub-linear expectation is introduced by Peng  \cite{PengG-Expectation06,peng2008a,peng2009survey} in a general function space by relaxing   the linear property of the classical linear expectation
to  the sub-additivity
and positive homogeneity (cf. Definition~\ref{def1.1} below).
   The sub-linear expectation does not depend on the probability measure,   provides
a  very flexible framework to model distribution uncertainty    problems  and produces many interesting properties different from those of the linear expectations. Under Peng's framework, many limit theorems have been being gradually established recently, including the central limit theorem and weak law of large numbers (cf. Peng \cite{peng2008a,peng2010}), the small derivation and Chung's law of the iterated logarithm (cf. Zhang \cite{Zhang Donsker}), the strong law of large numbers (cf. Chen \cite{chen2016strong}, Chen et al  \cite{Chen Z 2013}, Hu \cite{Hu C 2018}, Zhang \cite{Zhang Rosenthal}, Zhang and Lin \cite{ZhangLin}), and   the law of the iterated logarithm (cf.  Chen \cite{Chen2014LIL}, Zhang \cite{Zhang Exponential}).
For the convergence of the infinite series $\sum_{n=1}^{\infty}X_n$,  Xu and Zhang  \cite{XuZhang2018}  gave sufficient conditions of the almost sure convergence  for   independent random variables under the sub-linear expectation via a three-series theorem, recently. In this paper, we will consider the necessity of these conditions and the equivalence of the almost sure convergence, the convergence in capacity and the convergence in distribution. In the classical probability space, the Levy maximal inequalities are basic to the study of the almost sure behavior of sums of independent random variables and a key to show that the convergence in probability of $\sum_{n=1}^{\infty} X_n$ implies its almost sure convergence. We will establish  Levy type inequalities under the sub-linear expectation. For showing that the convergence in distribution of  $\sum_{n=1}^{\infty} X_n$ implies  its convergence in probability, the characteristic function is a basic tool. But, under the sub-linear expectation, there is no such tools.  We will find a new way to show a similar implication under the sub-linear expectations basing on a Komlogorov type maximal inequality.

As for the central limit theorem, it is well-known that the finite variances and mean zeros are sufficient and necessary for $\frac{\sum_{k=1}^n X_k}{\sqrt{n}}$ to converge in distribution to a normal random variable if $\{X_n;n\ge 1\}$ is a sequence of independent and identically distributed random variables on a classical probability space $(\Omega, \mathcal{F}, \pr)$. Under the sub-linear expectation, Peng \cite{peng2008a, peng2010} proved the cental limit theorem under  the finite $(2+\alpha)$-th moment. By applying a moment inequality and the truncation method, Zhang \cite{Zhang Exponential} and Lin and Zhang \cite{LinZhang2017} showed that the moment condition can be weakened to the finite second moment. A nature question is whether the finite second moment is necessary. In this paper, by applying the maximal inequalities, we will obtain the sufficient and necessary conditions for the central limit theorem.

In the remainder of the section, we state some natation.  In the next section, we will establish the maximal inequalities for random variables under the sub-linear expectation. The results on the convergence of the infinite series of random variables will given in Section \ref{Sect Convergece}. The sufficient and necessary conditions for the  central limit theorem are given in Section \ref{Sect CLT}.

We use the framework and notations of Peng  \cite{peng2008a}. Let  $(\Omega,\mathcal F)$
 be a given measurable space  and let $\mathscr{H}$ be a linear space of real functions
defined on $(\Omega,\mathcal F)$ such that if $X_1,\ldots, X_n \in \mathscr{H}$  then $\varphi(X_1,\ldots,X_n)\in \mathscr{H}$ for each
$\varphi\in C_{l,Lip}(\mathbb R_n)$,  where $C_{l,Lip}(\mathbb R_n)$ denotes the linear space of local Lipschitz
functions $\varphi$ satisfying
\begin{eqnarray*} & |\varphi(\bm x) - \varphi(\bm y)| \le  C(1 + |\bm x|^m + |\bm y|^m)|\bm x- \bm y|, \;\; \forall \bm x, \bm y \in \mathbb R_n,&\\
& \text {for some }  C > 0, m \in \mathbb  N \text{ depending on } \varphi. &
\end{eqnarray*}
$\mathscr{H}$ is considered as a space of ``random variables''. In this case we denote $X\in \mathscr{H}$. In the paper, we also denote $C_{b,Lip}(\mathbb R_n)$ the space of bounded  Lipschitz
functions, $C_b(\mathbb R_n)$ the space of bounded continuous functions, and $C_b^{1}(\mathbb R_n)$ the space of bounded continuous functions with bounded continuous derivations on $\mathbb R_n$.

\begin{definition}\label{def1.1} A  sub-linear expectation $\Sbep$ on $\mathscr{H}$  is a function $\Sbep: \mathscr{H}\to \overline{\mathbb R}$ satisfying the following properties: for all $X, Y \in \mathscr H$, we have
\begin{description}
  \item[\rm (a)]  Monotonicity: If $X \ge  Y$ then $\Sbep [X]\ge \Sbep [Y]$;
\item[\rm (b)] Constant preserving: $\Sbep [c] = c$;
\item[\rm (c)] Sub-additivity: $\Sbep[X+Y]\le \Sbep [X] +\Sbep [Y ]$ whenever $\Sbep [X] +\Sbep [Y ]$ is not of the form $+\infty-\infty$ or $-\infty+\infty$;
\item[\rm (d)] Positive homogeneity: $\Sbep [\lambda X] = \lambda \Sbep  [X]$, $\lambda\ge 0$.
 \end{description}
 Here $\overline{\mathbb R}=[-\infty, \infty]$. The triple $(\Omega, \mathscr{H}, \Sbep)$ is called a sub-linear expectation space. Give a sub-linear expectation $\Sbep $, let us denote the conjugate expectation $\cSbep$of $\Sbep$ by
$$ \cSbep[X]:=-\Sbep[-X], \;\; \forall X\in \mathscr{H}. $$
\end{definition}
From the definition, it is easily shown that    $\cSbep[X]\le \Sbep[X]$, $\Sbep[X+c]= \Sbep[X]+c$ and $\Sbep[X-Y]\ge \Sbep[X]-\Sbep[Y]$ for all
$X, Y\in \mathscr{H}$ with $\Sbep[Y]$ being finite. Further, if $\Sbep[|X|]$ is finite, then $\cSbep[X]$ and $\Sbep[X]$ are both finite.
\begin{definition}\label{def1.2}
\begin{description}
  \item[ \rm (i)] ({\em Identical distribution}) Let $\bm X_1$ and $\bm X_2$ be two $n$-dimensional random vectors defined
respectively in sub-linear expectation spaces $(\Omega_1, \mathscr{H}_1, \Sbep_1)$
  and $(\Omega_2, \mathscr{H}_2, \Sbep_2)$. They are called identically distributed, denoted by $\bm X_1\overset{d}= \bm X_2$  if
$$ \Sbep_1[\varphi(\bm X_1)]=\Sbep_2[\varphi(\bm X_2)], \;\; \forall \varphi\in C_{b,Lip}(\mathbb R_n), $$
where $C_{b,Lip}(\mathbb R_n)$ is the space of bounded Lipschitz functions.
\item[\rm (ii)] ({\em Independence})   In a sub-linear expectation space  $(\Omega, \mathscr{H}, \Sbep)$, a random vector $\bm Y =
(Y_1, \ldots, Y_n)$, $Y_i \in \mathscr{H}$ is said to be independent to another random vector $\bm X =
(X_1, \ldots, X_m)$ , $X_i \in \mathscr{H}$ under $\Sbep$  if
$$ \Sbep [\varphi(\bm X, \bm Y )] = \Sbep \big[\Sbep[\varphi(\bm x, \bm Y )]\big|_{\bm x=\bm X}\big], \;\; \forall   \varphi\in C_{b,Lip}(\mathbb R_m \times \mathbb R_n). $$

 Random variables $\{X_n; n\ge 1\}$
 are said to be independent, if  $X_{i+1}$ is independent to $(X_1,\ldots, X_i)$ for each $i\ge 1$.
 \end{description}
\end{definition}
In Peng \cite{peng2008a,peng2010,peng2010b}, the space of the test function $\varphi$ is $C_{l,Lip}(\mathbb R_n)$.
Here, the test function $\varphi$ in the definition is limit in the space of bounded Lipschitz functions.     When the considered random variables have finite moments of each order, i.e., $\Sbep[|X|^p]<\infty$ for each $p>0$, then the space of test functions $ C_{b,Lip}(\mathbb R_n)$     can be equivalently extended to  $C_{l,Lip}(\mathbb R_n)$.

A function $V:\mathcal{F}\to [0,1]$ is called a capacity if $V(\emptyset)=0$, $V(\Omega)=1$ and $V(A\cup B)\le V(A)+V(B)$ for all $A, B\in \mathcal{F}$.  Let $(\Omega, \mathscr{H}, \Sbep)$ be a sub-linear space.  We denote a pair $(\Capc,\cCapc)$ of capacities by
$$ \Capc(A):=\inf\{\Sbep[\xi]: I_A\le \xi, \xi\in\mathscr{H}\}, \;\; \cCapc(A):= 1-\Capc(A^c),\;\; \forall A\in \mathcal F, $$
where $A^c$  is the complement set of $A$.
Then
$$
 \Sbep[f]\le \Capc(A)\le \Sbep[g], \;\;\cSbep[f]\le \cCapc(A) \le \cSbep[g],\;\;
\text{ if } f\le I_A\le g, f,g \in \mathscr{H}.
$$
It is obvious that $\Capc$ is sub-additive, i.e.,  $\Capc(A\bigcup B)\le \Capc(A)+\Capc(B)$. But $\cCapc$ and $\cSbep$ are not. However, we have
$$
  \cCapc(A\bigcup B)\le \cCapc(A)+\Capc(B) \;\;\text{ and }\;\; \cSbep[X+Y]\le \cSbep[X]+\Sbep[Y]
$$
due to the fact that $\Capc(A^c\bigcap B^c)=\Capc(A^c\backslash B)\ge \Capc(A^c)-\Capc(B)$ and $\Sbep[-X-Y]\ge \Sbep[-X]-\Sbep[Y]$.
Further,
if $X$ is not in $\mathscr{H}$, we define $\Sbep[X]$ by
$$ \Sbep[X]=\inf\{\Sbep[Y]: X\le Y, \; Y\in \mathscr{H}\}. $$
Then $\Capc(A)=\Sbep[I_A]$.

\begin{definition}\label{def1.3}
(I) A function $V: \mathcal{F}\to  [0, 1]$ is called to be countably sub-additive if
\[
V\Big(\bigcup_{n=1}^{\infty}A_n\Big)\leq \sum_{n=1}^{\infty}V(A_n),\ \ \forall A_n\in \mathcal{F}.
\]
(II) A function $V: \mathcal{F}\to  [0, 1]$ is called to be  continuous  if it satisfies:\\
(i) Continuity from below: $V(A_n)\uparrow V(A)$ if $A_n\uparrow A$, where $A_n, A\in \mathcal{F}$.\\
(ii) Continuity from above: $V(A_n)\downarrow V(A)$ if $A_n \downarrow A$, where $A_n, A\in \mathcal{F}$.
\end{definition}
It is easily seen that a continuous capacity is countably sub-additive.

  \section{Maximal inequalities}\label{Sect Inequality}
  \setcounter{equation}{0}

  In this section, we establish several inequalities on the maximal sums.  The first one is   the Levy maximal inequality.

\begin{lemma}\label{LevyIneq} Let $X_1,\cdots, X_n$ be independent random variables in  a sub-linear expectation space $(\Omega, \mathscr{H}, \Sbep)$, $S_k=\sum_{i=1}^k X_i$, and $0<\alpha<1$ be a real number. If there exist real constants $\beta_{n, k}$  such that
$$ \Capc\left(S_k-S_n\ge \beta_{n,k}+\epsilon\right)\le \alpha, \text{ for all } \epsilon>0 \text{ and } k=1,\cdots ,n, $$
then
\begin{equation}\label{eqLIQ1}   (1-\alpha) \mathbb{V}\left(\max_{k\le n}(S_k -\beta_{n,k})> x+\epsilon\right)\le   \mathbb{V}\left(S_n>x\right), \text{ for all }x>0, \epsilon>0.
\end{equation}
If there exist real constants $\beta_{n, k}$  such that
$$ \Capc\left(|S_k-S_n|\ge \beta_{n,k}+\epsilon\right)\le \alpha, \text{ for all } \epsilon>0 \text{ and } k=1,\cdots ,n, $$
then
\begin{equation}\label{eqLIQ2}   (1-\alpha) \mathbb{V}\left(\max_{k\le n}(|S_k| -\beta_{n,k})> x+\epsilon\right)\le   \mathbb{V}\left(|S_n|>x\right), \text{ for all }x>0, \epsilon>0.
\end{equation}
\end{lemma}
{\bf Proof.} We only give the proof of (\ref{eqLIQ1}) since the proof of (\ref{eqLIQ2}) is similar.
Let $g_{\epsilon}(x)$ be a function with
\begin{equation}\label{eqproofLIQ.1} g_{\epsilon} \in C_b^1(\mathbb R) \; \text{ and } I_{\{x\ge \epsilon \}}\leq g_{\epsilon}(x)\leq I_{\{x\ge \epsilon/2\}}  \text{ for all } x,
\end{equation}
 where $0<\epsilon<1/2$, $ C_b^1(\mathbb R)$ is the space of bounded continuous function having bounded continuous derivations.   Denote $Z_k=g_{\epsilon}\left(S_k-\beta_{n,k}-x\right)$,   $Z_0=0$ and $\eta_k=\prod_{i=1}^k(1-Z_i)$. Then  $S_n-S_m$ is independent to $(Z_1,\ldots,Z_m)$, and
\begin{align*}
&(1-\alpha) I\{\max_{k\le n}(S_k -\beta_{n,k})> x+\epsilon\}\\
= & (1-\alpha)\left[1-\prod_{k=1}^n I\{ S_k-\beta_{n,k}- x\le \epsilon\}\right] \\
 \le & (1-\alpha)\left[1-\eta_n\right]=(1-\alpha)\left[\sum_{m=1}^n  \eta_{m-1}  Z_m\right]\\
 =& \sum_{m=1}^n  \eta_{m-1}  Z_m I\{S_m-S_n<\beta_{n,m}+\epsilon/2\}\\
  & +\sum_{m=1}^n  \eta_{m-1}Z_m\left[1-\alpha-I\{S_m-S_n<\beta_{n,m}+\epsilon/2\}\right]\\
 \le & \sum_{m=1}^n  \eta_{m-1}Z_m I\{S_n>x\}
   +\sum_{m=1}^n  \eta_{m-1}Z_m\left[-\alpha+I\{S_m-S_n\ge \beta_{n,m}+\epsilon/2\}\right]\\
  =&  I\{ S_n >x\}
   +\sum_{m=1}^n  \eta_{m-1}Z_m\left[-\alpha+g_{\epsilon/2}\left(  S_m-S_n-\beta_{n,m} \right) \right],
\end{align*}
where the second inequality above is due to the fact that on the event $\{Z_m\ne 0\}$ and $\{S_m-S_n< \beta_{n,m}+\epsilon/2\}$ we have $S_n\ge S_m-(S_m-S_n)>x$.  Note
$$ \Sbep\left[g_{\epsilon/2}\left(S_m-S_n-\beta_{n,m}\right) \right]
\le \mathbb{V}\left(S_m-S_n\ge \beta_{n,m}+\epsilon/4\right)\le \alpha. $$
By the independence,
\begin{align*}
& \Sbep\left[\eta_{m-1}Z_m\left[-\alpha+g_{\epsilon/2}\left(S_m-S_n-\beta_{n,m}\right) \right]\right]\\
=&\Sbep\left[\eta_{m-1}Z_m\left\{-\alpha+\Sbep\left[g_{\epsilon/2}\left(S_m-S_n-\beta_{n,m}\right) \right]\right\}\right]\le 0.
\end{align*}
By the sub-additivity of $\Sbep$, it follows that
\begin{align*}
 & (1-\alpha) \mathbb{V}\left(\max_{k\le n}(S_k -\beta_{n,k})> x+\epsilon\right)\\
 \le  & \mathbb{V}\left(S_n>x\right)
  +\sum_{m=1}^n \Sbep\left[ \eta_{m-1}Z_m\left[-\alpha+g_{\epsilon/2}\left(  S_m-S_n-\beta_{n,m} \right) \right]\right]
 \\
 \le &  \mathbb{V}\left(S_n>x\right).
 \end{align*}
 The proof is completed. $\Box$

 The second lemma is on  the Kolmogorov type inequality.

\begin{lemma}\label{KolIneq} Let $X_1,\cdots, X_n$ be independent random variables in  a sub-linear expectation space $(\Omega, \mathscr{H}, \Sbep)$. Let $S_k=\sum_{i=1}^k X_i$.
\begin{description}
  \item[\rm (i) ]  Suppose  $|X_k|\le c$,  $k=1,\cdots, n$. Then
\begin{equation}\label{eqKIQ1}   \mathbb{V}\left(\max_{k\le n}|S_k|> x \right)\ge 1 -\frac{(x+c)^2+2x\sum_{k=1}^n \big\{\big(\Sbep[X_k]\big)^++\big(\Sbep[-X_k]\big)^+\big\}}{\sum_{k=1}^n \Sbep[X_k^2]},
\end{equation}
for all $x>0$.
  \item[\rm (ii)] Suppose $ X_k \le c$,  $\Sbep[X_k]\ge 0$, $k=1,\cdots, n$.   Then
\begin{equation}\label{eqKIQ2.1}   \mathbb{V}\left(\max_{k\le n} S_k>x \right)\ge 1 -\frac{ x+c }{\sum_{k=1}^n \Sbep[X_k]}\; \text{ for all }x>0.
\end{equation}
\end{description}

\end{lemma}
{\bf Proof.} (i) Let $g_{\epsilon}$ be defined as in (\ref{eqproofLIQ.1}).  Denote $Z_k=g_{\epsilon}\left(|S_k|-x\right)$, $Z_0=0$, $\eta_k=\prod_{i=1}^k(1-Z_i)$.
Then $I\{|S_k|\ge x+\epsilon\}\le Z_k\le I\{|S_k|> x\}$. Also, $|S_{k-1}|< x+\epsilon$ and $|S_k|< |S_{k-1}|+|X_k|\le x+\epsilon+c$  on the event $\{\eta_{k-1}\ne 0\}$. So
\begin{align*}
S_{k-1}^2 \eta_{k-1}  + 2S_{k-1}X_k \eta_{k-1} +X_k^2 \eta_{k-1}
=& S_k^2 \eta_k  +S_k^2 \eta_{k-1}  Z_k   \\
\le &S_k^2 \eta_k  +(x+\epsilon+c)^2\left[\eta_{k-1} -\eta_k \right].
\end{align*}
Taking the summation over $k$ yields
\begin{align*}
&\left(\sum_{k=1}^n \Sbep[X_k^2]\right)\eta_n  +\sum_{k=1}^n\left(X_k^2-\Sbep[X_k^2]\right) \eta_{k-1}
\le  \sum_{k=1}^nX_k^2 \eta_{k-1}  \\
\le & S_n^2 \eta_n  +(x+\epsilon+c)^2\left[1-\eta_n \right]-2\sum_{k=1}^n S_{k-1}X_k \eta_{k-1} \\
\le &(x+\epsilon)^2 \eta_n  +(x+\epsilon+c)^2\left[1-\eta_n \right]-2\sum_{k=1}^n S_{k-1}X_k \eta_{k-1}  \\
\le & (x+\epsilon+c)^2 -2\sum_{k=1}^n S_{k-1}X_k \eta_{k-1} .
 \end{align*}
Write $B_n^2=\sum_{k=1}^n \Sbep[X_k^2]$.  It follows that
\begin{align*} & 1-\frac{ (x+\epsilon+c)^2 }{B_n^2}+\frac{ \sum_{k=1}^n\left(X_k^2-\Sbep[X_k^2]\right) \eta_{k-1}  }{B_n^2}\\
\le & 1-\eta_n +\frac{2}{B_n^2}\sum_{k=1}^n \left[X_k S_{k-1}^-\eta_{k-1} -X_kS_{k-1}^+\eta_{k-1} \right].
 \end{align*}
 Note $$\Sbep[X_k S_{k-1}^-\eta_{k-1} ]=\Sbep[X_k]\Sbep[ S_{k-1}^-\eta_{k-1} ]\le (x+\epsilon) \big(\Sbep[X_k]\big)^+,$$
 $$\Sbep[-X_k S_{k-1}^+\eta_{k-1} ]=\Sbep[-X_k]\Sbep[ S_{k-1}^+\eta_{k-1} ]\le (x+\epsilon) \big(\Sbep[-X_k]\big)^+$$
  and
 \begin{align}\label{eqpppoofKIQ1}
 & \Sbep\left[\sum_{k=1}^n\left(X_k^2-\Sbep[X_k^2]\right) \eta_{k-1}    \right] \nonumber\\
 = & \Sbep\left[\Sbep\left[\sum_{k=1}^n\left(X_k^2-\Sbep[X_k^2]\right) \eta_{k-1}   \Big|X_1,\cdots,X_{n-1}  \right]\right] \nonumber\\
 = & \Sbep\left[\sum_{k=1}^{n-1}\left(X_k^2-\Sbep[X_k^2]\right) \eta_{k-1}  +\eta_{n-1}\Sbep[X_n^2-\Sbep[X_n^2] ]  \right]\nonumber\\
 =&   \Sbep\left[\sum_{k=1}^{n-1}\left(X_k^2-\Sbep[X_k^2]\right) \eta_{k-1}     \right]=\cdots =0.
 \end{align}
 It follows that
\begin{align*}
&  1-\frac{ (x+\epsilon+c)^2 }{B_n^2}-\frac{ 2(x+\epsilon)\sum_{k=1}^n \big\{\big(\Sbep[X_k]\big)^++\big(\Sbep[-X_k]\big)^+\big\} }{B_n^2} \\
& \;\;  \le \Sbep\left[ 1-\eta_n\right]\le \Capc\left(\max_{k\le n} |S_k|> x\right).
\end{align*}
 By letting $\epsilon\to 0$, we obtain (\ref{eqKIQ1}). The proof of (i) is completed.

(ii)   Redefine $Z_k$ and $\eta_k$ by  $Z_k=g_{\epsilon}\left(S_k-x\right)$, $Z_0=0$, $\eta_k=\prod_{i=1}^k(1-Z_i)$.
Then $I\{S_k\ge x+\epsilon\}\le Z_k\le I\{S_k> x\}$. Also,  $ S_{k-1}< x+\epsilon$ and $ S_k =  S_{k-1}  + X_k < x+\epsilon+c$  on the event $\{\eta_{k-1}\ne 0\}$. So
$$
S_{k-1}  \eta_{k-1}  +  X_k  \eta_{k-1}
=  S_k  \eta_k  +S_k  \eta_{k-1} Z_k
\le S_k  \eta_k  +(x+\epsilon+c) \eta_{k-1} Z_k.
$$
Taking the summation over $k$ yields
\begin{align*}
&\left(\sum_{k=1}^n \Sbep[X_k ]\right)\eta_n  +\sum_{k=1}^n\left(X_k -\Sbep[X_k ]\right) \eta_{k-1} \\
\le & \sum_{k=1}^nX_k  \eta_{k-1}
\le   S_n  \eta_n  +(x+\epsilon+c) \left[1-\eta_n \right] \\
\le &(x+\epsilon)  \eta_n  +(x+\epsilon+c) \left[1-\eta_n \right]
\le   (x+\epsilon+c) .
 \end{align*}
Write $e_n =\sum_{k=1}^n \Sbep[X_k]$.  It follows that
\begin{align*}   1-\frac{ (x+\epsilon+c)  }{e_n}+\frac{ \sum_{k=1}^n\left(X_k-\Sbep[X_k]\right) \eta_{k-1} }{e_n}
\le   1-\eta_n.
 \end{align*}
 Note
 \begin{align*}
  \Sbep\left[\sum_{k=1}^n\left(X_k -\Sbep[X_k ]\right) \eta_{k-1}    \right]
 =    \Sbep\left[\sum_{k=1}^{n-1}\left(X_k -\Sbep[X_k ]\right) \eta_{k-1}     \right]=\cdots =0,
 \end{align*}
 similar to (\ref{eqpppoofKIQ1}). It follows that
 $$ 1-\frac{  x+\epsilon+c   }{e_n } \le \Sbep\left[ 1-\eta_n \right]\le \Capc\left(\max_{k\le n} S_k > x\right). $$
 By letting $\epsilon\to 0$, we obtain (\ref{eqKIQ1}). The proof is completed. $\Box$

The following lemma on the bounds of the capacities via moments will be used in the paper.
\begin{lemma}[ \cite{Zhang Exponential}]\label{moment_v}
Let $X_1, X_2, \ldots, X_n$ be  independent
random variables in $(\Omega, \mathscr{H}, \Capc)$. If $\Sbep[X_{k}] \leq 0$, $k=1,\ldots, n$, then  there exists a constant $C>0$ such that
\begin{equation*}
\Capc(S_n\geq x)\leq C\frac{\sum_{k=1}^{n}\Sbep[X_k^2]}{x^2} \;\text{ for all } \; \forall x>0.
\end{equation*}
\end{lemma}

  \section{The convergence of infinite series}\label{Sect Convergece}
  \setcounter{equation}{0}

Our  results on the convergence of the series $\sum_{n=1}^{\infty}$ are stated as three theorems. The first one gives the equivalency between the almost sure convergence and the convergence in capacity.
\begin{theorem}\label{th1} Let $\{X_n;n\geq1\}$ be a sequence of independent random variables in  a sub-linear expectation space $(\Omega, \mathscr{H}, \Sbep)$, $S_n=\sum_{k=1}^n X_k$, and $S$ be a random variable in the measurable space $(\Omega, \mathcal{F})$.
\begin{description}
  \item[(i)] If $\Capc$ is countably sub-additive, and
 \begin{equation}\label{eqth1.1} \Capc\left(|S_n-S|\ge \epsilon\right)\to 0 \text{ as } n\to \infty  \text{ for all } \epsilon>0,
 \end{equation}
  then
 \begin{equation}\label{eqth1.2} \Capc\left(\left\{\omega: \lim_{n\to \infty} S_n(\omega)\ne S(\omega)\right\}\right)=0.
 \end{equation}
 When (\ref{eqth1.2}) holds, we call that $\sum_{n=1}^{\infty}X_n$ is almost surely convergent in capacity, and when (\ref{eqth1.1}) holds, we call that $\sum_{n=1}^{\infty}X_n$ is  convergent in capacity.
  \item[(ii)] If $\Capc$ is continuous, then  (\ref{eqth1.2}) implies (\ref{eqth1.1}).
\end{description}
\end{theorem}

The second theorem gives the equivalency between the convergence in capacity and the convergence in distribution.
\begin{theorem}\label{th2} Let $\{X_n;n\geq1\}$ be a sequence of independent random variables in  a sub-linear expectation space $(\Omega, \mathscr{H}, \Sbep)$, $S_n=\sum_{k=1}^n X_k$.
\begin{description}
  \item[(i)] If there   is a random variable $S$ in the measurable space $(\Omega, \mathcal{F})$ such that
 \begin{equation}\label{eqth2.1} \Capc\left(|S_n-S|\ge \epsilon\right)\to 0 \text{ as } n\to \infty  \text{ for all } \epsilon>0,
 \end{equation}
 and $S$ is tight under $\Sbep$, i.e., $\Sbep\left[I_{\{|S|\le x\}^c}\right]=\Capc(|S|> x)\to 0$ as $x\to \infty$,
  then
 \begin{equation}\label{eqth2.2} \Sbep\left[\phi(S_n)\right]\to \Sbep\left[\phi(S)\right],\;\; \phi\in C_b(\mathbb{R}),
 \end{equation}
 where $C_b(\mathbb R)$ is the space of   bounded continuous functions on $\mathbb R$. When (\ref{eqth2.2}) holds, we call that $\sum_{n=1}^{\infty}X_n$ is  convergent in distribution.
  \item[(ii)] Suppose that there is a sub-linear space $(\widetilde{\Omega}, \widetilde{\mathscr{H}}, \widetilde{\mathbb E})$ and a random variable $\widetilde{S}$ on it such that $\widetilde{S}$ is tight under $\widetilde{\mathbb E}$, i.e., $\widetilde{\Capc}(|\widetilde{S}|> x)\to 0$ as $x\to \infty$, and
       \begin{equation}\label{eqth2.3} \Sbep\left[\phi(S_n)\right]\to \widetilde{\mathbb E} \left[\phi(\widetilde{S})\right],\;\; \phi\in C_{b}(\mathbb{R}),
 \end{equation}
      then $S_n$ is a Cauchy sequence in capacity $\Capc$, namely
      \begin{equation}\label{eqth2.4} \Capc\left(|S_n-S_m|\ge \epsilon \right)\to 0 \text{ as } n,m\to \infty \text{ for all } \epsilon>0.
      \end{equation}
      Furthermore, if $\Capc$ is countably sub-additive, then on the measurable space $(\Omega, \mathcal F)$ there is a random variable $S$ which is tight under $\Sbep$, such that (\ref{eqth1.1}) and (\ref{eqth1.2}) hold.
\end{description}
\end{theorem}

Recently, Xu and Zhang \cite{XuZhang2018} gave sufficient conditions for $\sum_{n=1}^{\infty} X_n$ to be convergent almost surely in capacity via three series theorem. The third  theorem of us gives   the sufficient and necessary conditions for $S_n$ to be a Cauchy sequence in capacity.
For any random variable $X$  and constant $c$, we denote  $X^c=(-c)\vee(X\wedge c)$.

\begin{theorem}\label{th4}  Let $\{X_n;n\geq1\}$ be a sequence of independent random variables in $(\Omega, \mathscr{H}, \Sbep)$, $S_n=\sum_{k=1}^n X_k$.    Then $S_n$ will be a Cauchy sequence in capacity $\Capc$ if the following three conditions hold for some $c>0$.
\begin{description}
  \item[\rm (S1) ] $\sum\limits_{n=1}^{\infty}\Capc(|X_n|>c)<\infty$,
  \item[\rm (S2) ] $\sum\limits_{n=1}^{\infty}\Sbep[X_n^c]$ and  $ \sum\limits_{n=1}^{\infty}\Sbep[-X_n^c]$ are both convergent,
  \item[\rm (S3) ]  $\sum\limits_{n=1}^{\infty}\Sbep\left[ \big(X_n^c-\Sbep[X_n^c]\big)^2\right] <\infty$ or/and $\sum\limits_{n=1}^{\infty}\Sbep\left[ \big(X_n^c+\Sbep[-X_n^c]\big)^2\right] <\infty$.
\end{description}
Conversely, if $S_n$  is  a Cauchy sequence in capacity $\Capc$, then (S1),(S2) and (S3) will hold for all $c>0$.
\end{theorem}

From Theorem \ref{th4}, we have the following three series theorem  on the sufficient and necessary conditions for the almost sure convergence of $\sum_{n=1}^{\infty} X_n$.
\begin{corollary}\label{three series} Let $\{X_n;n\geq1\}$ be a sequence of independent random variables in $(\Omega, \mathscr{H}, \Sbep)$. Suppose that  $\Capc$ is countably sub-additive. Then $\sum_{n=1}^{\infty}X_n$ will converge almost surely in capacity if the three conditions (S1),(S2) and (S3) in Theorem \ref{th4}  hold for some $c>0$.
Conversely, if  $\Capc$ is continuous  and $\sum_{n=1}^{\infty}X_n$ is  convergent almost surely in capacity, then (i),(ii) and (iii) will hold for all $c>0$.
\end{corollary}
The sufficiency of (S1), (S2) and (S3) is proved by Xu and Zhang \cite{XuZhang2018}, and also follows from  Theorem \ref{th4} and the second part of conclusion  of Theorem \ref{th2} (ii). The necessity follows from Theorem \ref{th4} and  Theorem \ref{th1} (ii).

  \bigskip

 The prove Theorems \ref{th1} and \ref{th2}. We need some more lemmas.
The first lemma is a version of Theorem 9 of Peng \cite{peng2010b}.
\begin{lemma}\label{lem1} Let $\{\bm Y_n; n\ge 1\}$ be a sequence of $d$-dimensional random variables in  a sub-linear expectation space $(\Omega, \mathscr{H}, \Sbep)$. Suppose that $\bm Y_n$ is asymptotically tight, i.e.,
$$  \limsup_{n\to\infty}\Sbep\left[I_{\{\bm Y_n\|\le x\}^c}\right]=\limsup_{n\to\infty} \Capc\left(\|\bm Y_n\|> x\right)\to 0 \; \text{ as } x\to \infty. $$
Then for any subsequence $\{\bm Y_{n_k}\}$  of $\{\bm Y_n\}$, there exist   further a subsequence $\{\bm Y_{n_{k^{\prime}}}\}$ of $\{\bm Y_{n_k}\}$ and a sub-linear expectation space $(\overline{\Omega}, \overline{\mathscr{H}}, \overline{\mathbb E})$ with a $d$-dimensional random variable $\bm Y$ on it such that
$$ \Sbep\left[\phi\left(\bm Y_{n_{k^{\prime}}}\right)\right]\to \overline{\mathbb E}\left[\phi(\bm Y)\right] \text{ for any } \phi\in C_b(\mathbb R^d) $$
and $\bm Y$ is tight under $\overline{\mathbb E}$.
\end{lemma}

{\bf Proof.} Let
$$ \mathbb E\left[\phi\right]=\limsup_{n\to \infty}\Sbep\left[\phi(\bm Y_n)\right], \;\; \phi\in C_b(\mathbb R^d). $$
Then $\mathbb E$ is a sub-linear expectation on the function space $C_b(\mathbb R^d)$ and is tight in sense that for any $\epsilon>0$, there is a compact set $K=\{\bm x:\|\bm x\|\le M\}$  for which
$\mathbb E\left[I_{K^c}\right]<\epsilon$. With the same argument as in the proof of Theorem 9 of Peng  \cite{peng2010b}, there is a countable subset $\{\varphi_j\}$ of $C_b(\mathbb R^d)$ such that for each $\phi\in C_b(\mathbb R^d)$ and any  $\epsilon>0$ one can find a  $\varphi_j$ satisfying
\begin{equation}\label{eqprooflem1.1}\mathbb E\left[|\phi-\varphi_j|\right]<\epsilon. \end{equation}
On the other hand, for each $\varphi_j$, the sequence $\Sbep\left[\varphi_j(\bm Y_n)\right]$ is bounded and so there is a Cauchy subsequence. Note that the set $\{\varphi_j\}$ is countable. By the diagonal choice method, one can find a sequence $\{n_k\}\subset \{n\}$ such that $\Sbep\left[\varphi_j(\bm Y_{n_k})\right]$ is a Cauchy sequence for each $\varphi_j$. Now, we show that $\Sbep\left[\phi(\bm Y_{n_k})\right]$ is a Cauchy sequence for any $\phi\in C_b(\mathbb R^d)$.  For any $\epsilon>0$, choose a $\varphi_j$ such that (\ref{eqprooflem1.1}) holds. Then
\begin{align*}
 & \left|\Sbep\left[\phi(\bm Y_{n_k})\right]-\Sbep\left[\phi(\bm Y_{n_l})\right]\right| \\
 \le & \left|\Sbep\left[\varphi_j(\bm Y_{n_k})\right]-\Sbep\left[\varphi_j(\bm Y_{n_l})\right]\right| \\
 &+  \Sbep\left[\left|\phi(\bm Y_{n_k}) - \varphi_j(\bm Y_{n_k})\right|\right]  + \Sbep\left[\left|\phi(\bm Y_{n_l}) -\varphi_j(\bm Y_{n_l})\right|\right].
 \end{align*}
 Taking the limits yields
 $$
\limsup_{k,l\to\infty}  \left|\Sbep\left[\phi(\bm Y_{n_k})\right]-\Sbep\left[\phi(\bm Y_{n_l})\right]\right|
 \le    0+ 2\mathbb E\left[|\phi-\varphi_j|\right]<2\epsilon.
$$
Hence $\Sbep\left[\phi(\bm Y_{n_k})\right]$ is a Cauchy sequence for any $\phi\in C_b(\mathbb R^d)$, and then
\begin{equation}\label{eqprooflem1.2} \lim_{k\to \infty} \Sbep\left[\phi(\bm Y_{n_k})\right] \; \text{ exists and is finite for any } \phi\in C_b(\mathbb R^d).
\end{equation}
Now, let $\overline{\Omega}=\mathbb R^d$, $\overline{\mathscr{H}}=C_{l,lip}(\mathbb R^d)$. Define
$$ \overline{\mathbb E}\left[\varphi\right]=\limsup_{k\to \infty} \Sbep\left[\varphi(\bm Y_{n_k})\right], \;\; \varphi\in C_{l,lip}(\mathbb R^d). $$
Then $(\overline{\Omega}, \overline{\mathscr{H}},\overline{\mathbb E})$ is a sub-linear expectation space. Define the random variable $\bm Y$ by $\bm Y(\bm x)=\bm x$, $\bm x\in \overline{\Omega}$.
From (\ref{eqprooflem1.2}) it follows that
$$ \lim_{k\to \infty} \Sbep\left[\phi(\bm Y_{n_k})\right]=\overline{\mathbb E}\left[\phi(\bm Y)\right] \text{ for any } \phi\in C_b(\mathbb R^d).
$$
The proof is completed. $\Box$

 \begin{lemma}\label{lem2}  Let $X$ and $Y$ be random variables in  a sub-linear expectation space $(\Omega, \mathscr{H}, \Sbep)$. Suppose that $Y$  and $X$ are independent ($Y$ is independent to $X$, or $X$ is independent to $Y$), and $X$ is tight, i.e. $\Capc(|X|\ge x)\to 0$ as $x> \infty$. If $X+Y\overset{d}=X$, then $\Capc(|Y|\ge \epsilon)=0$ for all $\epsilon>0$.
 \end{lemma}

 {\bf Proof.} Without loss of generality, we assume that $Y$ is independent to $X$. We can find a sub-linear expectation space $(\Omega^{\prime}, \mathscr{H}^{\prime}, \Sbep^{\prime})$ on which there are independent random variables  $X_1, Y_1,Y_2, \cdots, Y_n,\cdots$ such that $X_1\overset{d}=X$, $Y_i\overset{d}=Y$, $i=1,2,\cdots,$. Without loss of generality, assume $(\Omega^{\prime}, \mathscr{H}^{\prime}, \Sbep^{\prime})=(\Omega, \mathscr{H}, \Sbep)$. Let  $S_k=\sum_{j=1}^k Y_k$. Then $X_1+S_k\overset{d}=X$. So,
  \begin{align}\label{eqprooflem2.1}   \max_{k\le n} \Capc(|S_k|>x_0)\le  & \max_{k\le n} \Capc(|X_1+S_k|>x_0/2)+ \Capc(|X_1|>x_0/2) \nonumber
  \\
  \le &  \Sbep\left[g_{1/2}\left(\frac{|X_1+S_k|}{x_0}\right)\right]+ \Sbep\left[g_{1/2}\left(\frac{|X_1|}{x_0}\right)\right]\nonumber\\
  =& 2 \Sbep\left[g_{1/2}\left(\frac{|X|}{x_0}\right)\right]\le 2\Capc(|X|\ge  x_0/4)<1/4
  \end{align}
 for $x_0$ large enough,  where $g_{\epsilon}$ is defined as in (\ref{eqproofLIQ.1}).
 By Lemma \ref{LevyIneq},
\begin{equation}\label{eqprooflem2.2}  \Capc(\max_{k\le n}|S_k|>2x_0+\epsilon)\le \frac{4}{3} \max_n \Capc(|S_n|>x_0) \le \frac{4}{3} \cdot 2 \Capc(|X|\ge x_0/4)<\frac{1}{3} .
\end{equation}
 It follows that for any $\epsilon>0$,
 $$   \Capc(\max_{k\le n}|Y_k|>4x_0+2\epsilon) <\frac{1}{3}. $$
 Let $Z_k=g_{\epsilon} (|Y_k|-4x_0-2\epsilon)$, where $g_{\epsilon}$ is defined as in (\ref{eqproofLIQ.1}). Denote $q=\Capc(|Y_1|>4x_0+3\epsilon)$.  Then $Z_1,Z_2\cdots, Z_n$ are independent and identically distributed with $\{|Y_k|>4x_0+3\epsilon\}\le Z_k\le \{|Y_k|>4x_0+2\epsilon\}$ and  $\Sbep[Z_1]\ge \Capc(|Y_1|>4x_0+3\epsilon)=q$. Then by Lemma \ref{KolIneq} (ii),
 \begin{equation} \label{eqprooflem2.3} \frac{1}{3}> \Capc\left(\sum_{k=1}^n Z_k\ge 1\right)
 \ge 1-\frac{1+1}{\sum_{k=1}^n \Sbep[Z_k]}\ge 1-\frac{2}{nq}.
\end{equation}
 The above inequality holds for all $n$, which is impossible unless $q=0$. So we conclude that
 $$\Capc(|Y_1|>4x_0+\epsilon) = 0 \;\; \text{ for any } \epsilon>0. $$
 Now, let $\widetilde{Y}_k=(-5x_0)\vee Y_k \wedge (5x_0)$, $\widetilde{S}_k=\sum_{i=1}^k \widetilde{Y}_i$. Then $\widetilde{Y}_1, \cdots, \widetilde{Y}_n$ are independent and identically distributed  bounded random  variables, $\Capc(\widetilde{Y}_k\ne Y_k)=0$ and $\Capc(\widetilde{S}_k\ne S_k)=0 $.  If $\Sbep[\widetilde{Y}_1]>0$, then by Lemma \ref{KolIneq} (ii) again,
$$
\Capc(\max_{k\le n}S_k\ge 3x_0)=\Capc(\max_{k\le n}\widetilde{S}_k\ge 3x_0)\ge    1 -\frac{3x_0+5x_0}{n \Sbep[\widetilde{Y}_1]}, $$
which contradicts to (\ref{eqprooflem2.2}) when $n>  12 x_0/\Sbep[\widetilde{Y}_1]$. Hence, $\Sbep[\widetilde{Y}_1]\le 0$. Similarly, $\Sbep[-\widetilde{Y}_1]\le 0$. We conclude that $\Sbep[\widetilde{Y}_1]=\Sbep[-\widetilde{Y}_1]= 0$. Now, if $ \Sbep[\widetilde{Y}_1^2]\ne 0$, then by Lemma \ref{KolIneq} (i) we have
$$
\Capc(\max_{k\le n}|S_k|\ge 3x_0)\ge    1 -\frac{(3x_0+5x_0)^2}{n \Sbep[\widetilde{Y}_1^2]}, $$
which contradicts to (\ref{eqprooflem2.2}) when $n> 96x_0^2/\Sbep[\widetilde{Y}_1^2]$. We conclude that $\Sbep[\widetilde{Y}_1^2]=0$.

Finally, for any $\epsilon>0$ ($\epsilon<5x_0$),
$$ \Capc\left(|Y|\ge \epsilon\right)\le   \frac{\Sbep[Y^2\wedge(5x_0)^2]}{\epsilon^2} = \frac{\Sbep[\widetilde{Y}_1^2]}{\epsilon^2}=0. $$
The proof is completed. $\Box$

\bigskip

{\bf Proof of Theorem \ref{th1}.} (i) Let $\epsilon_k=1/2^k$, $\delta_k=1/4^k$. By (\ref{eqth1.1}), there exits a sequence $n_1<n_2<\cdots<n_k\to \infty$, such that
\begin{equation}\label{eqproofth1.1} \max_{n\ge n_k} \Capc\left(|S_n-S|\ge \epsilon_k\right)<\delta_k.
\end{equation}
By the countably sub-additivity of $\Capc$, we have
\begin{align*}
&\Capc\left(\limsup_{k\to \infty}|S_{n_k}-S|>0\right)
\le
 \Capc\left(\bigcap_{m=1}^{\infty} \bigcup_{k=m}^{\infty} \{|S_{n_k}-S|\ge \epsilon_k\}\right)\\
\le & \sum_{k=m}^{\infty}\Capc\left( |S_{n_k}-S|\ge \epsilon_k \right)\le \sum_{k=m}^{\infty}\delta_k \to 0 \text{ as } m\to \infty.
\end{align*}
By (\ref{eqproofth1.1}), $\max_{n\ge n_k} \Capc\left(|S_n-S_{n_{k+1}}|\ge 2\epsilon_k\right)<2\delta_k<1/2$. Apply the Levy inequality (\ref{eqLIQ2}) yields
\begin{equation}\label{eqproofth1.2} \Capc\left(\max_{n_k\le n\le n_{k+1}}|S_n-S_{n_k}|> 5\epsilon_k\right)\le 2  \Capc\left(|S_{n_{k+1}}-S_{n_k}|>  2\epsilon_k\right)<4 \delta_k.
\end{equation}
By the countably sub-additivity of $\Capc$ again,
\begin{align*}
&\Capc\left(\limsup_{k\to \infty}\max_{n_k\le n\le n_{k+1}}|S_n-S_{n_k}|>0\right)\\
\le &
\Capc\left(\bigcap_{m=1}^{\infty} \bigcup_{k=m}^{\infty} \{\max_{n_k\le n\le n_{k+1}}|S_n-S_{n_k}|\ge 5\epsilon_k\}\right)\\
\le & \sum_{k=m}^{\infty}\Capc\left( \max_{n_k\le n\le n_{k+1}}|S_n-S_{n_k}|\ge 5\epsilon_k \right)\le 4\sum_{k=m}^{\infty}\delta_k \to 0 \text{ as } m\to \infty.
\end{align*}
It follows that
\begin{align*}
& \Capc\left(\limsup_{n\to \infty}|S_n-S|>0\right)\\
\le & \Capc\left(\limsup_{k\to \infty}|S_{n_k}-S|>0\right)+\Capc\left(\limsup_{k\to \infty}\max_{n_k\le n\le n_{k+1}}|S_n-S_{n_k}|>0\right)=0.
\end{align*}
(\ref{eqth1.2}) follows.

(ii) From (\ref{eqth1.2}) and the continuity of $\Capc$, it follows that for any $\epsilon>0$,
\begin{align*}
0\ge \Capc\left(\bigcap_{n=1}^{\infty}\bigcup_{m=n}^{\infty}\{|S_m-S|\ge \epsilon\}\right)= & \lim_{n\to \infty}\Capc\left(\bigcup_{m=n}^{\infty}\{|S_m-S|\ge \epsilon\}\right) \\
\ge  & \limsup_{n\to \infty}\Capc\left( |S_n-S|\ge \epsilon \right).
\end{align*}
(\ref{eqth1.1}) follows.  The proof is completed. $\Box$

{\bf Proof of Theorem \ref{th2}.} (i) We first show that (\ref{eqth2.2}) holds for any bounded uniformly continuous function $\phi$. For any $\epsilon>0$, there is a $\delta>0$ such that $|\phi(x)-\phi(y)|<\epsilon$ when $|x-y|<\delta$. It follows that
$$ \left|\Sbep\left[\phi(S_n)\right]-\Sbep\left[\phi(S)\right]\right|\le \epsilon+ 2\sup_x|\phi(x)|\Capc\left(|S_n-S|\ge \delta\right). $$
By letting $n\to \infty$ and the arbitrariness of $\epsilon>0$, we obtain (\ref{eqth2.2}). Now, suppose that $\phi$ is a bounded continuous function. Then for any $N>1$, $\phi((-N)\vee x\wedge N)$ is a bounded uniformly continuous function. Hence
$$ \lim_{n\to \infty}\Sbep\left[\phi((-N)\vee S_n\wedge N)\right]=\Sbep\left[\phi((-N)\vee S\wedge N)\right]. $$
On the other hand,
$$\left|\Sbep\left[\phi((-N)\vee S\wedge N)\right]-\Sbep\left[\phi(  S )\right]\right|\le 2\sup_x\big|\phi(x)\big|\Capc\left(|S|>N\right)\to 0 \text{ as } N\to \infty, $$
and
\begin{align*}
&\limsup_{n\to \infty} \left|\Sbep\left[\phi((-N)\vee S_n\wedge N)\right]-\Sbep\left[\phi(  S_n )\right]\right|\\
\le & 2\sup_x|\phi(x)|\limsup_{n\to \infty} \Capc\left(|S_n|\ge N\right)\le 2\sup_x\big|\phi(x)\big|\limsup_{n\to \infty} \Sbep\left[g_1\left(\frac{|S_n|}{N}\right)\right]\\
=& 2\sup_x|\phi(x)|  \Sbep\left[g_1\left(\frac{|S|}{N}\right)\right]\le 2\sup_x\big|\phi(x)\big|\limsup_{n\to \infty} \Capc\left(|S|\ge N/2\right)\to 0 \text{ as } N\to \infty,
\end{align*}
where $g_{\epsilon}$ is defined as in (\ref{eqproofLIQ.1}).  Hence, (\ref{eqth2.2}) holds for a bounded continuous function $\phi$.

(ii) Note
$$ \Capc\left(|S_n-S_m|\ge 2x\right)\le \Capc\left(|S_n|\ge x\right)+\Capc\left(|S_m|\ge x\right). $$
It follows that
\begin{align*}
&\limsup_{m\ge n\to \infty}\Capc\left(|S_n-S_m|\ge 2x\right)\le 2\limsup_{n\to \infty} \Capc\left(|S_n|\ge x\right) \\
\le & 2\limsup_{n\to \infty} \Sbep\left[g_1\left(\frac{|S_n|}{x}\right)\right]=\widetilde{\mathbb E} \left[g_1\left(\frac{|\widetilde{S}|}{x}\right)\right]
\le 2\Capc\left(|\widetilde{S}|\ge x/2\right) \to 0 \text{ as } x\to \infty.
\end{align*}
Write $\bm Y_{n,m}=(S_n, S_m-S_n)$, then the sequence $\{\bm Y_{n,m}; m\ge n\}$ is asymptotically tight, i.e.,
$$\limsup_{m\ge n\to \infty} \Capc\left(\|\bm Y_{n,m}\|\ge x\right)\to 0\; \text{ as } x \to \infty. $$
By Lemma \ref{lem1}, for any subsequence $(n_k,m_k)$ of $(n,m)$, there is  further a subsequence $(n_{k^{\prime}},m_{k^{\prime}})$ of $(n_k,m_k)$ and a sub-linear expectation space $(\overline{\Omega}, \overline{\mathscr{H}}, \overline{\mathbb E})$ with a random vector $\bm Y=(Y_1,Y_2)$ such that
\begin{equation}\label{eqproofth2.3} \Sbep\left[\phi\left(\bm Y_{n_{k^{\prime}},m_{k^{\prime}}}\right)\right]\to \overline{\mathbb E}\left[\phi(\bm Y)\right],\;\;\phi\in C_b(\mathbb R^2).
\end{equation}
Note that $S_{m_{k^{\prime}}}- S_{n_{k^{\prime}}}$ is independent to $S_{n_{k^{\prime}}}$. By Lemma 4.4 of Zhang  \cite{Zhang Donsker}, $Y_2$ is independent to $Y_1$ under $ \overline{\mathbb E}$.
Let $\phi\in C_{b,Lip}(\mathbb R)$. By (\ref{eqproofth2.3}),
\begin{equation}\label{eqproofth2.4}\Sbep\left[\phi\left(S_{m_{k^{\prime}}}\right)\right]\to \overline{\mathbb E}\left[\phi(Y_1+Y_2)\right], \;\; \Sbep\left[\phi\left(S_{n_{k^{\prime}}}\right)\right]\to \overline{\mathbb E}\left[\phi(Y_1)\right]
\end{equation}
and
\begin{equation}\label{eqproofth2.5} \Sbep\left[\phi\left(S_{m_{k^{\prime}}}-S_{n_{k^{\prime}}}\right)\right]\to \overline{\mathbb E}\left[\phi(Y_2)\right].
\end{equation}
On the other hand, by (\ref{eqth2.3}),
\begin{equation}\label{eqproofth2.6}  \Sbep\left[\phi\left(S_{m_{k^{\prime}}}\right)\right] \to \widetilde{\mathbb E} \left[\phi(\widetilde{S})\right]
\; \text{ and } \Sbep\left[\phi\left(S_{n_{k^{\prime}}}\right)\right] \to \widetilde{\mathbb E} \left[\phi(\widetilde{S})\right].
\end{equation}
Combing (\ref{eqproofth2.4}) and (\ref{eqproofth2.6}) yields
$$ \overline{\mathbb E}\left[\phi(Y_1+Y_2)\right] =\overline{\mathbb E}\left[\phi(Y_1)\right]=\widetilde{\mathbb E} \left[\phi(\widetilde{S})\right],\;\;\phi\in C_{b,Lip}(\mathbb R). $$
Hence, by Lemma \ref{lem2}, we obtain $\overline{\mathbb V}(|Y_2|\ge \epsilon)=0$ for all $\epsilon>0$. By choosing $\phi\in C_{b,Lip}(\mathbb R)$ such that $I_{|x|\ge \epsilon}\le \phi(x)\le I_{|x|\ge \epsilon/2}$ in
(\ref{eqproofth2.5}), we have
$$
\limsup_{k^{\prime}\to \infty} \Capc\left(\left|S_{m_{k^{\prime}}}-S_{n_{k^{\prime}}}\right|\ge \epsilon\right)\le \overline{\mathbb V}(|Y_2|\ge \epsilon/2)=0.
$$
So, we conclude that for any subsequence $(n_k,m_k)$ of $(n,m)$, there is a further a subsequence $(n_{k^{\prime}},m_{k^{\prime}})$ of $(n_k,m_k)$ such that
$$ \Capc\left(\left|S_{m_{k^{\prime}}}-S_{n_{k^{\prime}}}\right|\ge \epsilon\right)\to 0 \text{ for all }\epsilon>0. $$
Hence (\ref{eqth2.4}) is proved.

Next, suppose that $\Capc$ is countably sub-additive. Let $\epsilon_k=1/2^k$, $\delta_k=1/3^k$. By (\ref{eqth2.4}), there is a sequence $n_1<n_2<\cdots<n_k<\cdots$ such that
$$ \Capc\left(|S_{n_{k+1}}-S_{n_k}|\ge \epsilon_k\right)\le \delta_k.  $$
Let $A=\{\omega: \sum_{k=1}^{\infty}|S_{n_{k+1}}-S_{n_k}|<\infty\}$. Then
\begin{align*} \Capc\left(A^c\right)\le & \Capc\left(\sum_{k=K}^{\infty}|S_{n_{k+1}}-S_{n_k}|\ge \sum_{k=K}^{\infty}\epsilon_k\right)\\
\le & \sum_{k=K}^{\infty} \Capc\left(|S_{n_{k+1}}-S_{n_k}|\ge \epsilon_k\right)\le \sum_{k=K}^{\infty}\delta_k \to 0 \text{ as } K\to \infty.
\end{align*}
Define $S=\lim_{k\to \infty} S_{n_k}$ on $A$, and $S=0$ on $A^c$. Then
\begin{align*} \Capc\left(|S-S_{n_k}|\ge 1/2^{k-1}\right)\le &\Capc(A^c)+ \Capc\left(A, \sum_{i=k}^{\infty}|S_{n_{i+1}}-S_{n_i}|\ge \sum_{i=k}^{\infty}\epsilon_i\right)\\
\le & \sum_{i=k}^{\infty} \Capc\left(|S_{n_{i+1}}-S_{n_i}|\ge \epsilon_i\right)\le \sum_{i=k}^{\infty}\delta_i \to 0 \text{ as } k\to \infty.
\end{align*}
On the other hand, by (\ref{eqth2.4}),
$$\Capc\left(|S_n-S_{n_k}|\ge \epsilon\right)\to 0 \text{ as } n, n_k\to \infty. $$
Hence
$$\Capc\left(|S_n-S|\ge \epsilon\right)\le  \Capc\left(|S_n-S_{n_k}|\ge \epsilon/2\right)+\Capc\left(|S-S_{n_k}|\ge \epsilon/2\right)\to 0. $$
(\ref{eqth1.1}) is proved. Further,
\begin{align*}
\Capc\left(|S|\ge 2M\right)\le & \limsup_n \Capc\left(|S_n|\ge M\right)+\limsup_n \Capc\left(|S_n-S|\ge M\right)\\
\le  & \widetilde{\mathbb V}\left(|\widetilde{S}|\ge M/2\right)\to 0
\text{ as } M\to \infty.
\end{align*}
So, $S$ is tight.
Finally, (\ref{eqth1.2}) follows from Theorem \ref{th1}. $\Box$

\bigskip

For showing Theorem \ref{th4}, we need a more lemma.
\begin{lemma}\label{lem3} Let $\{X_n; n\ge 1\}$  be a sequence of independent random variables in  a sub-linear expectation space $(\Omega, \mathscr{H}, \Sbep)$ with $ |X_k| \le c$, $\Sbep[X_k]\ge 0$ and $\Sbep[-X_k]\ge 0$, $k=1,2,\cdots$. Let $S_k=\sum_{i=1}^k X_i$. Suppose
\begin{equation}\label{eqlem3.1}   \lim_{x\to \infty}\lim_{n\to \infty} \mathbb{V}\left(\max_{k\le n} |S_k|>x \right)<1.
\end{equation}
Then  $\sum_{n=1}^{\infty}\Sbep[X_n]$, $\sum_{n=1}^{\infty}\Sbep[-X_n]$ and $\sum_{n=1}^{\infty}\Sbep[X_n^2]$ are convergent.
\end{lemma}

{\bf Proof.}     By (\ref{eqlem3.1}), there exist $0<\beta<1$, $x_0>0$ and $n_0$, such that
$$\mathbb{V}\left(\max_{k\le n} |S_k|>x \right) <\beta, \;\;\text{for all } x\ge x_0, \; n\ge n_0. $$
By (\ref{eqKIQ2.1}),
$$ \sum_{k=1}^n \Sbep[X_k]\le \frac{x+c}{1-\beta},\;\;  \text{for all } x\ge x_0, \; n\ge n_0. $$
So $\sum_{k=1}^{\infty}  \Sbep[X_k]$ is convergent. Similarly, $\sum_{k=1}^{\infty}  \Sbep[-X_k]$ is convergent.

Now, by (\ref{eqKIQ1}),
\begin{align*} \sum_{k=1}^n \Sbep[X_k^2] \le &  \frac{(x+c)^2+2x\sum_{k=1}^n \big\{\big(\Sbep[X_k]\big)^++\big(\Sbep[-X_k]\big)^+\big\}}{1-\beta} \\
 \le & \frac{(x+c)^2+2x\sum_{k=1}^{\infty}   \big\{ \Sbep[X_k] + \Sbep[-X_k] \big\}}{1-\beta},\;\;  \text{for all } x\ge x_0, \; n\ge n_0.
\end{align*}
So $\sum_{n=1}^{\infty}  \Sbep[X_n^2]$ is convergent. The proof   is completed.
$\Box$

\bigskip
{\bf Proof of Theorem \ref{th4}.}  (i) By Lemma \ref{moment_v} and the condition (S3),
\begin{align*}
& \Capc\left(S_n-S_m-\sum_{k=m+1}^n \Sbep[X_k]\ge \epsilon\right)\\
\le & C\frac{\sum_{k=m+1}^n\Sbep\left[(X_k-\Sbep[X_k])^2\right]}{\epsilon^2}\to 0 \text{ as } n\ge m\to \infty.
\end{align*}
The convergence of $\sum_{n=1}^{\infty} \Sbep[X_n] $ implies $\sum_{k=m+1}^n \Sbep[X_k]\to 0$. It follows that
$$ \lim_{n\ge m\to \infty}\Capc\left(S_n-S_m \ge \epsilon\right)=0\; \text{ for all }\epsilon>0. $$
On the other hand, note $\Sbep[X_k]+\Sbep[-X_k]\ge 0$. The condition (S2) implies $\sum_{n=1}^{\infty}\left(\Sbep[X_k]+\Sbep[-X_k]\right)<\infty$, and then $\sum_{n=1}^{\infty}\left(\Sbep[X_k]+\Sbep[-X_k]\right)^2<\infty$. Hence, by  the condition (S3) and the fact that $\Sbep\left[(-X_k-\Sbep[-X_k])^2\right]\le \Sbep\left[(X_k-\Sbep[X_k])^2\right]+(\Sbep[X_k]+\Sbep[-X_k])^2$,
$$\sum_{n=1}^{\infty}\Sbep\left[(-X_n-\Sbep[-X_n])^2\right]<\infty. $$
By considering $-X_n$ instead of $X_n$, we have
$$ \lim_{n\ge m\to \infty}\Capc\left(-S_n+S_m \ge \epsilon\right)=0\; \text{ for all }\epsilon>0. $$
It follows that (\ref{eqth2.4}) holds, i.e., $S_n$ is a Cauchy sequence in capacity $\Capc$.

(ii) Suppose that $S_n$ is a Cauchy sequence in capacity $\Capc$. Similar to (\ref{eqproofth1.2}), by applying the Levy inequality (\ref{eqLIQ2}) we have
\begin{equation}\label{eqproofth4.1} \lim_{n\ge m\to \infty} \Capc\left(\max_{m\le k\le n} |S_k-S_m|> \epsilon \right)=0 \text{ for all } \epsilon>0.
\end{equation}
Then
\begin{equation}\label{eqproofth4.2}\lim_{n\ge m\to \infty} \Capc\left(\max_{m\le k\le n} |X_k|\ge c \right)=0 \text{ for all } c>0.
\end{equation}
Write $v_k=\Capc\left( |X_k|\ge 2c \right)$. Similar to (\ref{eqprooflem2.3}), we have for $m_0$ large enough and all $n\ge m\ge m_0$,
\begin{align*}
\frac{1}{3}>\Capc\left(\max_{m\le k\le n} |X_k|\ge c \right)\ge 1-\frac{2}{ \sum_{k=m+1}^n v_k}.
\end{align*}
It follows that $\sum_{k=1}^{\infty}  v_k<\infty$. The condition (S1) is satisfied for all $c>0$.

Next, we consider (S3). Write $X_n^c=(-c)\vee X_n\wedge c$ and $S_n^c=\sum_{k=1}^n X_k^c$. Note on the event $\{\max_{m\le k\le n} |X_k|< c\}$, $X_k^c=X_k$, $k=m+1,\cdots, n$. By (\ref{eqproofth4.1}) and (\ref{eqproofth4.2}),
\begin{equation}\label{eqproofth4.3} \lim_{n\ge m\to \infty} \Capc\left(\max_{m\le k\le n} |S_k^c-S_m^c|> \epsilon \right)=0 \text{ for all } \epsilon>0.
\end{equation}
 Let $Y_1, Y_1^{\prime}, Y_2, Y_2^{\prime},\cdots, Y_n, Y_n^{\prime}, \cdots $ be independent random variables under the sub-linear expectation $\Sbep$ with $Y_k\overset{d}=Y_k^{\prime}\overset{d}= X_k^c$, $k=1,2,\cdots$. Then
 $$\{Y_{m+1},\cdots, Y_n \} \overset{d}=\{Y_{m+1}^{\prime},\cdots, Y_n^{\prime}\} \overset{d}=\{X_{m+1}^c,\cdots, X_n^c\}. $$
  Let $T_k=\sum_{i=1}^k Y_i$ and $T_k^{\prime}=\sum_{i=1}^k Y_i^{\prime}$. By (\ref{eqproofth4.3}),
 \begin{align}\label{eqproofth4.4} & \lim_{n\ge m\to \infty} \Capc\left(\max_{m\le k\le n} |T_k-T_m|> \epsilon \right)\nonumber \\
 = &\lim_{n\ge m\to \infty} \Capc\left(\max_{m\le k\le n} |T_k^{\prime}-T_m^{\prime}|> \epsilon \right)=0 \text{ for all } \epsilon>0.
\end{align}
 Write $\widetilde{Y}_n=Y_n-Y_n^{\prime}$ and $\widetilde{T}_n=\sum_{k=1}^n \widetilde{Y}_k$. Then $\{\widetilde{Y}_n;n\ge 1\}$ is a sequence of independent random variables with $\Capc(|\widetilde{Y}_n|> 3c)=0$. Without loss of generality, we can assume $|\widetilde{Y}_n|\le 3c$ for otherwise we can replace $\widetilde{Y}_n$ by $(-3c)\vee \widetilde{Y}_n \wedge(3c)$. By (\ref{eqproofth4.4}),
$$ \lim_{n\ge m\to \infty} \mathbb{V}\left(\max_{m\le k\le n} |\widetilde{T}_k-\widetilde{T}_m|>2\epsilon \right)=0 \text{ for all } \epsilon>0.$$
Note $\Sbep[-\widetilde{Y}_k]=\Sbep[\widetilde{Y}_k]=(\Sbep[X_k^c]+\Sbep[-X_k^c])/2\ge 0$. By Lemma \ref{lem3},
$$ \sum_{n=1}^{\infty} \left(\Sbep[X_n^c]+\Sbep[-X_n^c]\right)\;\text{ and } \sum_{n=1}^{\infty} \Sbep[\widetilde{Y}_n^2] \text{ are convergent}. $$
Note
\begin{align*}
\Sbep\left[\widetilde{Y}_n^2|Y_n\right]\ge & \big(Y_n-\Sbep[Y_n]\big)^2+\Sbep\left[\big(Y_n^{\prime}-\Sbep[Y_n^{\prime}]\big)^2\right]\\
 &+2\big(Y_n-\Sbep[Y_n]\big)^-\cSbep\left[Y_n^{\prime}-\Sbep[Y_n^{\prime}] \right].
\end{align*}
So
\begin{align*}
\Sbep\left[\widetilde{Y}_n^2\right]\ge  & 2\Sbep\big[ \big(X_n^c-\Sbep[X_n^c]\big)^2\big]-2 \{\Sbep[X_n^c]+\Sbep[-X_n^c]\} \Sbep\left[\big(X_n^c-\Sbep[X_n^c]\big)^-\right]
\\
\ge  & 2\Sbep\big[ \big(X_n^c-\Sbep[X_n^c]\big)^2\big]-2c \{\Sbep[X_n^c]+\Sbep[-X_n^c]\}.
\end{align*}
It follows that
\begin{equation} \label{eqproofth4.5}
\sum_{n=1}^{\infty}\Sbep\big[ \big(X_n^c-\Sbep[X_n^c]\big)^2\big]<\infty. \end{equation}
Since $\Sbep\big[ \big(-X_n^c-\Sbep[-X_n^c]\big)^2\big]\le \Sbep\big[ \big(X_n^c-\Sbep[X_n^c]\big)^2\big]+ \big(\Sbep[X_n^c+\Sbep[-X_n^c]\big)^2$, we also have
$$ \sum_{n=1}^{\infty}\Sbep\big[ \big(-X_n-\Sbep[-X_n]\big)^2\big]<\infty. $$
The condition (S3) is proved.

Finally, we consider (S2). For any $\epsilon>0$, when $m,n$ are large enough, $\sum_{k=m+1}^n \big(\Sbep[X_n^c]+\Sbep[-X_n^c]\big)<\epsilon$. By (\ref{eqproofth4.5}) and Lemma \ref{moment_v},
\begin{align*}
&\Capc\left(S_n^c-S_m^c-\sum_{k=m+1}^n \frac{\Sbep[X_k^c]-\Sbep[-X_k^c]}{2}\ge \epsilon\right)\\
= & \Capc\left(S_n^c-S_m^c-\sum_{k=m+1}^n \Sbep[X_k^c] \ge \epsilon-\sum_{k=m+1}^n \frac{\Sbep[-X_k^c]+\Sbep[X_k^c]}{2}\right)\\
\le & C \frac{\sum_{k=m+1}^n \Sbep\big[ \big(X_k^c-\Sbep[X_k^c]\big)^2\big]}{(\epsilon/2)^2}\to 0 \text{ as } n\ge m\to \infty.
\end{align*}
Similarly, by considering $-X_k^c$ instead of $X_k^c$ we have
\begin{align*}
 \Capc\left(-S_n^c+S_m^c-\sum_{k=m+1}^n \frac{\Sbep[-X_k^c]-\Sbep[X_k^c]}{2}\ge \epsilon\right) \to 0 \text{ as } n\ge m\to \infty.
\end{align*}
It follows that, for any $\epsilon>0$,
$$
\Capc\left(\left|S_n^c-S_m^c-\sum_{k=m+1}^n \frac{\Sbep[X_k^c]-\Sbep[-X_k^c]}{2}\right|\ge \epsilon\right) \to 0 \text{ as } n\ge m\to \infty,
$$
which, together with (\ref{eqproofth4.3}), implies
$$ \sum_{k=m+1}^n \frac{\Sbep[X_k^c]-\Sbep[-X_k^c]}{2} \to 0 \text{ as } n\ge m\to \infty. $$
Hence, $\sum_{n=1}^{\infty} \big(\Sbep[X_k^c]-\Sbep[-X_k^c]\big)$ is convergent. Note that $\sum_{n=1}^{\infty} \big(\Sbep[X_k^c]+\Sbep[-X_k^c]\big)$ is convergent. We conclude that both
$\sum_{n=1}^{\infty} \Sbep[X_k^c]$ and $\sum_{n=1}^{\infty}  \Sbep[-X_k^c]$ are convergent. The proof of (ii) is completed. $\Box$.

  \section{Central limit theorem}\label{Sect CLT}
  \setcounter{equation}{0}
In this section, we consider the sufficient and necessary conditions for the central limit theorem. We first recall the definition of G-normal random variables  which is introduced by Peng \cite{peng2008a, peng2010}.

\begin{definition}\label{def4.1} ({\em G-normal random variable})
For $0\le \underline{\sigma}^2\le \overline{\sigma}^2<\infty$, a random variable $\xi$ in a sub-linear expectation space $(\widetilde{\Omega}, \widetilde{\mathscr H}, \widetilde{\mathbb E})$   is called a normal $N\big(0, [\underline{\sigma}^2, \overline{\sigma}^2]\big)$ distributed  random variable   (written as $\xi \sim N\big(0, [\underline{\sigma}^2, \overline{\sigma}^2]\big)$  under $\widetilde{\mathbb E}$), if for any  $\varphi\in C_{l,Lip}(\mathbb R)$, the function $u(x,t)=\widetilde{\mathbb E}\left[\varphi\left(x+\sqrt{t} \xi\right)\right]$ ($x\in \mathbb R, t\ge 0$) is the unique viscosity solution of  the following heat equation:
  \begin{equation}\label{eqheatequation}\partial_t u -G\left( \partial_{xx}^2 u\right) =0, \;\; u(0,x)=\varphi(x),
  \end{equation}
where $G(\alpha)=\frac{1}{2}(\overline{\sigma}^2 \alpha^+ - \underline{\sigma}^2 \alpha^-)$.
\end{definition}
That $\xi$ is a normal distributed random variable is equivalent to that,   if $\xi^{\prime}$ is an independent copy of $\xi$ (i.e., $\xi^{\prime}$ is independent to $\xi$ and $\xi\overset{d}=\xi^{\prime})$, then
\begin{equation}\label{eqnormal} \widetilde{\mathbb E}\left[\varphi(\alpha \xi+\beta \xi^{\prime})\right]
=\widetilde{\mathbb E}\left[\varphi\big(\sqrt{\alpha^2+\beta^2}\xi\big)\right], \;\;
\forall \varphi\in C_{l,Lip}(\mathbb R) \text{ and } \forall \alpha,\beta\ge 0,
\end{equation}
(cf. Definition II.1.4 and Example II.1.13 of Peng \cite{peng2010}). We also write $\eta\overset{d}=  N\big(0, [\underline{\sigma}^2, \overline{\sigma}^2]\big)$ if $\eta\overset{d}=\xi$ (as defined in Definition \ref{def1.2} (i)) and $\xi \sim N\big(0, [\underline{\sigma}^2, \overline{\sigma}^2]\big)$ (as defined in Definition \ref{def4.1}). By definition, $\eta\overset{d}=\xi$ if and only if for any  $\varphi\in C_{b,Lip}(\mathbb R)$, the function $u(x,t)=\Sbep\left[\varphi\left(x+\sqrt{t} \eta\right)\right]$ ($x\in \mathbb R, t\ge 0$) is the unique viscosity solution of  the equation (\ref{eqheatequation}). In the sequel, without loss of generality, we assume that the sub-linear expectation spaces $(\widetilde{\Omega}, \widetilde{\mathscr H}, \widetilde{\mathbb E})$ and $(\Omega, \mathscr{H},\Sbep)$ are the same.

Let $\{X_n; n\ge 1\}$ be a sequence of independent and identically distributed random variables in a sub-linear expectation space $(\Omega, \mathscr{H},\Sbep)$, $S_n=\sum_{k=1}^nX_k$. Peng \cite{peng2008a, peng2010} proved that, if $\Sbep[X_1]=\Sbep[-X_1]=0$ and $\Sbep[|X_1|^{2+\alpha}]<\infty$ for some $\alpha>0$, then
 \begin{equation}\label{cltpeng}
    \lim_{n\to \infty} \Sbep\left[\varphi\left(\frac{S_n}{\sqrt{n}}\right)\right]=\Sbep\left[\varphi(\xi )\right], \forall \varphi\in C_b(\mathbb R),
    \end{equation}
    where $\xi\sim N\left(0,[\underline{\sigma}^2,\overline{\sigma}^2]\right)$, $\overline{\sigma}^2=\Sbep[X_1^2]$ and $\underline{\sigma}^2=\cSbep[X_1^2]$. Zhang \cite{Zhang Exponential} showed that
    $\Sbep[|X_1|^{2+\alpha}]<\infty$ can be weakened to $\Sbep[(X_1^2-c)^+]\to 0$ as $c\to\infty$ by applying the moment inequalities of sums of independent random variables and the truncation method. A nature question is whether   $\Sbep[X_1^2]<\infty$ and $ \Sbep[X_1]=\Sbep[-X_1]=0$ are   sufficient and necessary for (\ref{cltpeng}). The following theorem is our main result.

  \begin{theorem}\label{thclt} Let $\{X_n; n\ge 1\}$ be a sequence of independent and identically distributed random variables in a sub-linear expectation space $(\Omega, \mathscr{H},\Sbep)$, $S_n=\sum_{k=1}^nX_k$. Suppose that
  \begin{description}
    \item[\rm (i) ] $\lim_{c\to\infty} \Sbep[X_1^2\wedge c]$ is finite;
    \item[\rm (ii)]  $x^2\Capc\left(|X_1|\ge x\right)\to 0$ as $x\to \infty$;
    \item[\rm (iii)]  $\lim_{c\to \infty}\Sbep\left[(-c)\vee X_1\wedge c)\right]=\lim_{c\to \infty}\Sbep\left[(-c)\vee (- X_1)\wedge c)\right]=0$.
     \end{description}
    Write $\overline{\sigma}^2=\lim_{c\to\infty} \Sbep[X_1^2\wedge c]$ and  $\underline{\sigma}^2=\lim_{c\to\infty} \cSbep[X_1^2\wedge c]$. Then for any $\varphi\in C_b(\mathbb R)$,
    \begin{equation}\label{clt1}
    \lim_{n\to \infty} \Sbep\left[\varphi\left(\frac{S_n}{\sqrt{n}}\right)\right]=\Sbep\left[\varphi(\xi )\right],
    \end{equation}
    where $\xi\sim N\left(0,[\underline{\sigma}^2,\overline{\sigma}^2]\right)$.

    Conversely, if (\ref{clt1}) holds for any $\varphi\in C_b^1(\mathbb R)$ and a   random variable $\xi$ with  $x^2\Capc\left(|\xi|\ge x\right)\to 0$ as $x\to \infty$, then
    (i),(ii) and (iii) hold and  $\xi\overset{d}= N\left(0,[\underline{\sigma}^2,\overline{\sigma}^2]\right)$.
 \end{theorem}

Before prove the theorem, we give some remarks on the conditions. Note that $\Sbep[X_1^2\wedge c]$ and $\cSbep[X_1^2\wedge c]$ are non-decreasing in $c$. So, $\overline{\sigma}^2$ and $\underline{\sigma}^2$ are well-defined and nonnegative, and are finite if the condition (i) is satisfied. It is easily seen that, for $c_1>c_2>0$,
\begin{equation}\label{eqproofclt6}\left|\Sbep[X_1^{c_1}]-\Sbep[X_1^{c_2}]\right|\le \Sbep[(|X_1|\wedge c_1-c_2)^+]\le \frac{\overline{\sigma}^2}{c_2}.
\end{equation}
So, the condition (i) implies that $ \lim_{c\to \infty}\Sbep[X_1^{c}]$ and $ \lim_{c\to \infty}\Sbep[-X_1^{c}]$ exist  and are finite.

If $\Sbep$ is a continuous sub-linear expectation, i.e., $\Sbep[X_n]\nearrow \Sbep[X]$ whenever $0\le X_n\nearrow X$, and $\Sbep[X_n]\searrow 0$ whenever $X_n\searrow 0$, $\Sbep[X_n]<\infty$,
 then (i) is equivalent to $\Sbep[X_1^2]<\infty$, (iii) is equivalent to $\Sbep[X_1]=\Sbep[-X_1]=0$, and (ii) is automatically implied by $\Sbep[X_1^2]<\infty$. In general, the condition $\Sbep[X_1^2]<\infty$ and (i) with (ii) do not imply  each other. However, it is easily verified that, if $\Sbep[(X_1^2-c)^+]\to 0$ as $c\to \infty$, then (i) and (ii) are satisfied and (iii) is equivalent to $\Sbep[X_1]=\Sbep[-X_1]=0$.

 \bigskip
 To prove Theorem \ref{thclt}, we need a more lemma.

 \begin{lemma}\label{lem4.1} Let $X_{n1},\cdots X_{nn}$  be   independent   random variables in a sub-linear expectation space $(\Omega, \mathscr{H},\Sbep)$ with
 $$ \frac{1}{\sqrt{n}}\sum_{k=1}^n \left\{\left|\Sbep[X_{nk}]\right|+\left|\Sbep[-X_{nk}]\right|\right\}\to 0, $$
 $$ \frac{1}{n}\sum_{k=1}^n \left\{\big|\Sbep[X_{nk}^2]-\overline{\sigma}^2\big|+\big|\cSbep[X_{nk}^2]-\underline{\sigma}^2\big|\right\}\to 0 $$
 and
 $$ \frac{1}{n^{3/2}}\sum_{k=1}^n  \Sbep[|X_{nk}|^3]\to 0. $$
 Then
 $$\lim_{n\to \infty} \Sbep\left[\varphi\left(\frac{\sum_{k=1}^n X_{nk}}{\sqrt{n}}\right)\right]=\Sbep\left[\varphi(\xi )\right], \forall \varphi\in C_b(\mathbb R), $$
 where $\xi\sim N\left(0,[\underline{\sigma}^2,\overline{\sigma}^2]\right)$.
 \end{lemma}

 This lemma can be proved by refining the arguments of Li and Shi \cite{LiShi10} and can also follow from the Lindeberg central limit theorem \cite{Zhang Lindeberg}. We omit the proof here.

 \bigskip

 {\bf Proof of Theorem \ref{thclt}. } We first prove the sufficient part, i.e., (i),(ii) and (iii) $\implies$ (\ref{clt1}). Let $X_{nk}= (-\sqrt{n})\vee   X_k   \wedge \sqrt{n}$. Then
 for any $\epsilon>0$,
 $$\frac{1}{n^{3/2}} \sum_{k=1}^n \Sbep[|X_{nk}|^3]=\frac{1}{n^{1/2}}   \Sbep[|X_{n1}|^3]\le \epsilon \overline{\sigma}^2+n\Capc\left(|X_1|\ge \epsilon\sqrt{n}\right)\to 0
 $$
 as $n\to \infty$ and then $\epsilon\to 0$,
 by the condition (ii).    Also,
\begin{align*}
& \frac{1}{n}\sum_{k=1}^n \left\{\big|\Sbep[X_{nk}^2]-\overline{\sigma}^2\big|+\big|\cSbep[X_{nk}^2]-\underline{\sigma}^2\big|\right\}\\
=&  \big|\Sbep\left[X_1^2\wedge n\right]
 -\overline{\sigma}^2\big|+\big|\cSbep\left[X_1^2\wedge n\right]
 -\underline{\sigma}^2\big|\to 0,
 \end{align*}
 by (i).  Note by (ii) and (i),
 \begin{align*}
  &\frac{1}{\sqrt{n}} \sum_{k=1}^n \left|\Sbep[X_{nk}]\right|=\sqrt{n}  \left|\Sbep[X_{n1}]\right|\\
  = & \sqrt{n} \lim_{c\to \infty} \left|\Sbep[X_{n1}]-\Sbep\left[(-c\sqrt{n})\vee  X_1  \wedge (c\sqrt{n})\right]\right]\\
 \le &  \sqrt{n}\lim_{c\to \infty}  \Sbep\left[ \left(|X_1| \wedge (c\sqrt{n})-x\sqrt{n}\right)^+\right]+ \sqrt{n}\Sbep\left[ \left(|X_1| \wedge (x\sqrt{n})- \sqrt{n}\right)^+\right]\\
 \le & \frac{\overline{\sigma}^2}{x}+ x n\Capc\left(|X_1|\ge \sqrt{n}\right)\to 0 \; \text{ as } n\to \infty \text{ and then } x\to \infty,
\end{align*}
and similarly,
$$ \frac{1}{\sqrt{n}}\sum_{k=1}^n \left|\Sbep[-X_{nk}]\right|\to 0. $$
The conditions in Lemma \ref{lem4.1} are satisfied.  We obtain
$$\lim_{n\to \infty} \Sbep\left[\varphi\left(\frac{\sum_{k=1}^n X_{nk}}{\sqrt{n}}\right)\right]=\Sbep\left[\varphi(\xi )\right]. $$
It is obvious that
$$\Sbep\left[\left|\varphi\left(\frac{\sum_{k=1}^n X_{nk}}{\sqrt{n}}\right)-\varphi\left(\frac{S_n}{\sqrt{n}}\right)\right|\right]\le \sup_x|\varphi(x)|n\Capc\left(|X_1|\ge \sqrt{n}\right) \to 0. $$
(\ref{clt1}) is proved.

Now, we consider the necessary part. Letting $\varphi=g_{\epsilon}\big(|x|-t\big)$ yields
$$ \limsup_{n\to \infty}\Capc\left(\frac{|S_n|}{\sqrt{n}}\ge t+\epsilon\right)\le \Capc\left(|\xi|\ge t \right) \text{ for all } t>0, \epsilon>0. $$
So
$$ \limsup_{n\ge m\to \infty}\max_{m\le k,l\le n} \Capc\left(\frac{|S_k-S_l|}{\sqrt{n}}\ge 2t+\epsilon\right)\le 2\Capc\left(|\xi|\ge t \right) \text{ for all } t>0, \epsilon>0. $$
Choose $t_0$ such that $\Capc\left(|\xi |\ge t_0\right)<1/(32)$. Applying the Levy maximal inequality (\ref{eqLIQ2}) yields
\begin{equation}\label{eqproofclt1} \limsup_{n\ge m\to \infty}  \Capc\left(\frac{\max_{m\le k\le n}|S_k-S_m|}{\sqrt{n}}\ge 4t\right)< \frac{64}{31}\Capc\left(|\xi|\ge t \right) \text{ for all } t>t_0. \end{equation}
Hence
\begin{equation}\label{eqproofclt2} \limsup_{n\ge m\to \infty}  \Capc\left(\frac{\max_{m\le k\le n}|X_k|}{\sqrt{n}}\ge 8t\right)< \frac{64}{31}\Capc\left(|\xi|\ge t \right) \text{ for all } t>t_0.
\end{equation}
Let $t_1>t_0$ and $m_0$ such that
\begin{equation}\label{eqproofclt3}  \Capc\left(\frac{\max_{m\le k\le n}|S_k-S_m|}{\sqrt{n}}> 4t_1\right)< \frac{2}{31}  \text{ for all } m\ge m_0
\end{equation}
and
\begin{equation}\label{eqproofclt4}   \Capc\left(\frac{\max_{m\le k\le n}|X_k|}{\sqrt{n}}> 8t_1\right)< \frac{4}{31}  \text{ for all } m\ge m_0.
\end{equation}
Write $Y_{nk}=(-8t_1)\vee\left(\frac{X_k}{\sqrt{n}}\right)\wedge(8t_1)$.  Then by (\ref{eqproofclt3}) and (\ref{eqproofclt4}),
\begin{equation}\label{eqproofclt5}  \Capc\left( \max_{m\le k\le n}\big|\sum_{j=m+1}^k Y_{nj}\big| > 4t_1\right)< \frac{2}{31} +  \frac{4}{31}<\frac{1}{5} \text{ for all } m\ge m_0
\end{equation}
If $\Sbep[Y_{n1}]>0$, then by Lemma \ref{KolIneq} (ii),
$$\frac{1}{5}>1-\frac{4t_1+8t_1}{(n-m) \Sbep[Y_{n1}]}. $$
Hence  $(n-m) \big(\Sbep[Y_{n1}])^+\le 15t_1$. Similarly, $(n-m) \big(\Sbep[-Y_{n1}])^+\le 15t_1$. Hence, by Lemma \ref{KolIneq} (i), it follows that
$$ \frac{1}{5}> 1-\frac{ (4t_1+8t_1)^2+8t_1 \big\{(n-m) \big(\Sbep[Y_{n1}])^++(n-m) \big(\Sbep[-Y_{n1}])^+\big\}}{ (n-m) \Sbep[Y_{n1}^2]}. $$
We conclude that $(n-m)\Sbep[Y_{n1}^2]\le \frac{5}{4}(12^2+240) t_1^2$. Choose $m=n/2$ and let $n\to \infty$. We have
$$ \lim_{c\to \infty}\Sbep[X_1^2\wedge c]=\lim_{n\to \infty} n \Sbep[Y_{n1}^2]\le \frac{5}{2}(12^2+240) t_1^2. $$
(i) is proved. Note that (i) implies that $ \lim_{c\to \infty}\Sbep[X_1^{c}]$ exists and is finite. Then
$$ \lim_{c\to \infty}\Sbep[X_1^{c}]=\limsup_{n\to \infty}\sqrt{n}\Sbep[Y_{n1}] \le \limsup_{n\to \infty}\frac{30t_1}{\sqrt{n}}=0. $$
Similarly, $ \lim_{c\to \infty}\Sbep[-X_1^{c}]$ exists,  is finite and  not positive. Note $\Sbep[-X_1^{c}]+\Sbep[X_1^{c}]\ge 0$.
Hence (iii) follows.

Finally, we show (ii). For any given $0<\epsilon<1/2$,  by the condition $x^2\Capc(|\xi|\ge x)\to 0$,  one can choose $t_1>t_0$ such that
$\frac{64}{31}\Capc(|\xi|\ge t_1)\le \frac{\epsilon}{9^3t_1^2}<1/2$. Then by (\ref{eqproofclt2}), there is $m_0$ such that
$$ \Capc\left(\frac{\max_{m\le k\le n}|X_k|}{\sqrt{n}}\ge 8t_1\right)< \frac{\epsilon}{9^3t_1^2}, \; n\ge m\ge m_0. $$
Choose $Z_k=g_{\epsilon}\big(\frac{|X_k|}{8t_1\sqrt{n}}-1\big)$ such that $I\{|X_k|\ge 9t_1\sqrt{n}\}\le Z_k\le I\{|X_k|\ge 8t_1\sqrt{n}\}$. Let $q_n=\Capc\left(|X_1|\ge 9t_1\sqrt{n}\right)$.  Then
\begin{align*}
\Capc\left(\frac{\max_{m\le k\le n}|X_k|}{\sqrt{n}}\ge 8t_1\right)\ge & \Sbep\left[1-\prod_{k=m+1}^n(1-Z_k)\right]\\
=& 1-\prod_{k=m+1}^n(1-\Sbep[Z_k])\ge 1-e^{-(n-m)q_n}.
\end{align*}
It follows that
$$ n\Capc\left(|X_1|\ge 9t_1\sqrt{n}\right)\le 2(n-m)q_n< 2\times 2\times \frac{\epsilon}{9^3t_1^2}\text{ for } m=[n/2]\ge m_0. $$
Hence
$$ (9t_1\sqrt{n})^2 \Capc\left(|X_1|\ge 9t_1\sqrt{n}\right)< \frac{4 \epsilon}{9}, \;\; n\ge 2m_0. $$
When $x\ge 9t_1 \sqrt{2m_0}$, there is $n$ such that $9t_1\sqrt{n}\le x\le 9t_1\sqrt{n+1}$. Then
$$ x^2\Capc\left(|X_1|\ge x\right)\le (9t_1\sqrt{n+1})^2 \Capc\left(|X_1|\ge 9t_1\sqrt{n}\right)\le \frac{8 \epsilon}{9}. $$
It follows that $\limsup_{x\to \infty}x^2\Capc\left(|X_1|\ge x\right)<\epsilon$. (ii) is proved. The proof is now completed. $\Box$

\begin{remark} From the proof, we can find that
$$ \lim_{x\to \infty}\limsup_{n\to \infty}\Capc\left(\frac{|S_n|}{\sqrt{n}}\ge x\right)=0 $$
implies (i) and (ii). One may conjecture  that,
  \begin{description}
    \item[\rm  C1]  if (\ref{clt1}) holds for any $\varphi\in C_b^1(\mathbb R)$ and a   tight random variable $\xi$ (i.e.,   $\Capc\left(|\xi|\ge x\right)\to 0$ as $x\to \infty$), then (i), (ii) and (iii)  holds and $\xi\overset{d}= N\left(0,[\underline{\sigma}^2,\overline{\sigma}^2]\right)$.
  \end{description}
An equivalent conjecture is that,
\begin{description}
    \item[\rm  C2]   if  $\xi$ and $\xi^{\prime}$ are independent and identically distributed tight random variables,  and
\begin{equation}\label{eqnormal2} \Sbep\left[\varphi(\alpha \xi+\beta \xi^{\prime})\right]
=\Sbep\left[\varphi\big(\sqrt{\alpha^2+\beta^2}\xi\big)\right], \;\;
\forall \varphi\in C_b(\mathbb R) \text{ and } \forall \alpha,\beta\ge 0,
\end{equation}
then $\xi\overset{d}= N\big(0,[\underline{\sigma}^2,\overline{\sigma}^2]\big)$, where $\overline{\sigma}^2=\lim_{c\to \infty}\Sbep[\xi^2\wedge c]$ and $\underline{\sigma}^2=\lim_{c\to \infty}\cSbep[\xi^2\wedge c]$.
 \end{description}
  It should be noted that the conditions   (\ref{eqnormal}) and (\ref{eqnormal2}) are different. The condition (\ref{eqnormal}) implies that $\xi$ have finite moments  of each order, but non  information about the moments of $\xi$ is  hidden in (\ref{eqnormal2}). As Theorem \ref{thclt}, the conjecture C2 is  true when $x^2\Capc\left(|\xi|\ge x\right)\to 0$ as $x\to \infty$. In fact, let $X_1,X_2, \cdots, $ be independent random variables with $X_k\overset{d}=\xi$. Then by (\ref{eqnormal2}), $\frac{S_n}{\sqrt{n}}\overset{d}=\xi$. By the necessary part of Theorem \ref{thclt}, the conditions (i), (ii) and (iii) are satisfied. Then by the sufficient part of the theorem, $\xi\overset{d}=N\big(0,[\underline{\sigma}^2,\overline{\sigma}^2]\big)$.   We don't known whether
     conjectures C1 and C2 are true   without assuming any moment conditions. It is very possible that they are not true in general. But finding a counterexample is not an easy task.
\end{remark}


\bigskip

\begin{thebibliography}{99}

\bibitem{chen2016strong}
Chen, Z. J.: Strong laws of large numbers for sub-linear expectation. \emph{Sci. China Math.}, \textbf{59}(5), 945--954 (2016)

\bibitem{Chen2014LIL}
Chen, Z. J., Hu, F.: A law of the iterated logarithm under sublinear expectations. \emph{Journal of Financial Engineering}, \textbf{1}, No.2 (2014)


\bibitem{Chen Z 2013}
Chen, Z. J., Wu, P. Y., Li, B. M.: A strong law of large numbers for non-additive probabilities. \emph{Int. J. Approx. Reason.}, \textbf{54}(3), 365--377 (2013)

 \bibitem{LiShi10} Li, M. and Shi, Y.F.:   A general central limit theorem under sublinear expectations. \emph{Science in China Ser. A}, \textbf{53}(8), 1989-1994.

 \bibitem{LinZhang2017} Lin, Z Y.  and Zhang, L.X.:  Convergence to a self-normalized G-Brownian motion. \emph{ Probability, Uncertainty and Quantitative Risk}, \textbf{2}(4).   doi:10.1186/s41546-017-0013-8

 \bibitem{Hu C 2018} Hu, C. : Strong laws of large numbers for sublinear expectation under controlled 1st moment condition. \emph{Chinese Annals of Mathematics Series B}, \textbf{39} (5):791-804.


\bibitem{PengG-Expectation06}
Peng, S. G.: G-expectation, G-Brownian motion and related stochastic calculus of Ito type. In: Proceedings of the 2005 Abel Symposium. Springer, Berlin-Heidelberg, 541--567, 2007


\bibitem{peng2008a}
Peng, S. G.: A new central limit theorem under sublinear expectations. ArXiv:0803.2656v1 (2008)

\bibitem{peng2010}
Peng, S. G.: Nonlinear Expectations and Stochastic Calculus under Uncertainty, arXiv:1002.4546 [math.PR] (2010)

\bibitem{peng2010b}
Peng, S. G.: Tightness, weak compactness of nonlinear
expectations and application to CLT, 	arXiv:1006.2541 [math.PR] (2010)

\bibitem{peng2009survey}
Peng, S. G.: Survey on normal distributions, central limit theorem, Brownian motion and the related stochastic calculus under sublinear expectations. \emph{Sci. China Ser. A}, \textbf{52}(7), 1391--1411 (2009)



\bibitem{XuZhang2018}
Xu, J. P. and Zhang, L. X.: Three series theorem for independent random variables underSub-linear expectations with applications, \emph{Acta Mathematica Sinica, English Series}, Published online:  https://doi.org/10.1007/s10114-018-7508-9 (2018)

\bibitem{Zhang Donsker}
Zhang, L. X.: Donsker's invariance principle under the sub-linear expectation with an application to Chung's law of the iterated logarithm, \emph{Communications in Math. Stat.}, \textbf{3}(2), 187--214 (2015)

\bibitem{Zhang Exponential}
Zhang, L. X.: Exponential inequalities under the sub-linear expectations with applications to laws of the
iterated logarithm. \emph{Sci. China Math.}, \textbf{59}(12), 2503--2526 (2016)


\bibitem{Zhang Rosenthal}
Zhang, L. X.: Rosenthal's inequalities for independent and negatively dependent random variables under sub-linear expectations with applications. \emph{Sci. China Math.}, \textbf{59}(4), 751--768 (2016)

\bibitem{Zhang Lindeberg}
Zhang, L. X.: Lindeberg's central limit theorems for martingale like sequences under nonlinear expectations,  arXiv:1611.01619 [math.PR] (2016)

\bibitem{ZhangLin}
Zhang, L. X.   and Lin, J. H.:  Marcinkiewicz's strong law of large numbers for nonlinear expectations. \emph{Stat. Probab. Lett.}, \textbf{137}, 269--276 (2018)

\end{thebibliography}
 \end{document}